\documentclass{article}

\begin{document}

\title{K\"ahler-Einstein metrics and stability}
\author{Xiu-Xiong Chen, Simon Donaldson, Song Sun}
\maketitle


\newcommand{\bC}{{\bf C}}
\newcommand{\bP}{{\bf P}}
\newcommand{\Fut}{{\rm Fut}}
\newtheorem{defn}{Definition}
\newtheorem{thm}{Theorem}

\section{Introduction}
This is a note to announce  a solution of the well-known question: which  Fano manifolds admit  K\"ahler-Einstein metrics? The idea that the appropriate condition should be in terms of \lq\lq algebro-geometric stability'' was proposed by Yau about 20 years ago \cite{kn:Y}(partly by analogy with the  \lq\lq Kobayashi-Hitchin correspondence'' in the case of holomorphic bundles).
Over the years various different notions of stability have been discussed in the literature, both in the Fano/K\"ahler-Einstein case and in the more general situation of constant scalar curvature K\"ahler metrics on polarised manifolds. But the condition we end up using is essentially as proposed by Tian in \cite{kn:Ti1} and which we now recall (see also \cite{kn:LX} for the equivalence with {\it a priori} stronger definitions). 

\begin{defn} Let $X$ be an $n$-dimensional Fano manifold. A test-configuration for $X$ is a flat family $\pi:{\cal X}\rightarrow \bC$ embedded in $\bC\bP^{N}\times C$ for some $N$, invariant under a $\bC^{*}$ action on $\bC\bP^{N}\times \bC$ covering the standard action on $\bC$ such that 
\begin{itemize}
\item $\pi^{-1}(1)=X$ and the embedding $X\subset \bC\bP^{N}$ is defined by the complete linear system $\vert -r K_{X}\vert $ for some $r$;
\item The central fibre $X_{0}=\pi^{-1}(0) $ is a normal variety with log terminal singularities.
\end{itemize}
\end{defn}

A test configuration has a basic numerical invariant: the {\it Futaki invariant}.
This can be defined in various ways. One way is to consider the hyperplane bundle $L\rightarrow X_{0}$ (restricted from the hyperplane bundle on $\bC\bP^{N}$).
For integers $k\geq 1$ we have a vector space $H^{0}(X_{0}, L^{k})$ with a $\bC^{*}$-action. Let $d_{k}$ be the dimension of this vector space and $w_{k}$ be the total weight of the action. By general theory these are, for large $k$, given by polynomials in $k$ of degrees $n, n+1$ respectively. Thus 
   \begin{equation} \frac{w_{k}}{k d_{k}} = F_{0} + F_{1} k^{-1} + O(k^{-2})\end{equation}
   and the Futaki invariant $\Fut({\cal X})$ can be defined to be the co-efficient $F_{1}$.
   
   \begin{defn} $X$ is K-stable if for all non-trivial test configurations ${\cal X}$ we have $\Fut({\cal X}\geq 0$ and strict inequality holds if ${\cal X}$ is non-trivial.
\end{defn}

By \lq\lq non-trivial'' here we mean that the central fibre is not isomorphic to $X$. The condition defined above is  often called \lq\lq polystability'' in the literature.
 Our main result is
 \begin{thm}
 If a  Fano manifold $X$ is $K$-stable then it admits  a K\"ahler-Einstein metric.
\end{thm}
  
  Converse results, in different degrees of generality, have been proved by Tian \cite{kn:Ti1}, Stoppa \cite{kn:Stoppa} and Berman \cite{kn:Ber2}. The sharp form proved by Berman shows that in fact K-stability (as we have defined it) is equivalent to the existence of a K\"ahler-Einstein metric. 
 
   \section{Outline of proof}
   \subsection{Strategy}
   The strategy of proof follows that suggested in \cite{kn:sd2}. We fix some $\lambda>0$ such that the linear system $\vert -\lambda K_{X}\vert$ contains a smooth divisor $D$. It seems that in all known cases one could take $\lambda=1$ which was the situation discussed in \cite{kn:sd2}. Now we consider K\"ahler-Einstein metrics on $X$ with a cone singularity of cone angle $\beta$ along $D$, where $\beta\in (0,1]$ is a variable parameter. Of course when $\beta=1$ these are just smooth K\"ahler-Einstein metrics.  Such metrics with cone singularities were discussed in general terms some years ago by Tian in \cite{kn:Ti2}. More recently, following \cite{kn:sd2}, a detailed theory has been developed, both on the differential geometric side \cite{kn:sd2}, \cite{kn:Brendle}, \cite{kn:JMR}, \cite{kn:SW} and the algebro-geometric side \cite{kn:S}, \cite{kn:LS}, \cite{kn:OS}. 

 If ${\cal X}$ is a test configuration for $X$, as above then we can extend $D\subset X=\pi^{-1}(1)$ to a divisor in ${\cal X}$ and obtain a $\bC^{*}$-invariant divisor $D_{0}\subset X_{0}$. Then there is a modified Futaki invariant
$$  \Fut_{\beta}({\cal X})= \Fut({\cal X})-2\pi (1-\beta) F_{0}(D_{0})$$
where $F_{0}(D_{0})$ is the numerical  invariant of the $\bC^{*}$-action defined as in (1), but using the divisor $D_{0}$ in place of $X_{0}$.
A simple but fundamental point is that this is linear in $\beta$, so if we know that $  \Fut_{\beta}({\cal X})=0$ for some $\beta\leq 1$ and we know that $\Fut_{\beta'}({\cal X})\geq 0$ for some $\beta'<\beta$ then we can deduce that $\Fut({\cal X})\leq 0$ and so $X$ is not K-stable. (There is a small extra complication in the case when $X$ admits holomorphic automorphisms which we will return to below.)

The existence of a K\"ahler-Einstein metric with small cone angle is well-understood. When $\lambda=1$ this is proved by Berman in \cite{kn:Ber1}. A particularly simple case is when $\lambda=2$, for then we get a Ricci-flat metric with $\beta=1/2$ by applying Yau's theorem to the double cover branched over $D$. In general for $\lambda>2$ we can use the existence theorems of Brendle \cite{kn:Brendle} or Jeffres, Mazzeo, Rubinstein \cite{kn:JMR}. See the discussion in \cite{kn:LS}. Likewise the \lq\lq K-semistability''---the fact that for any test configuration we have $\Fut_{\beta'}({\cal X})\geq 0$---is well understood for small $\beta'$. \cite{kn:S}\cite{kn:OS}. Thus it suffices  to show that if $X$ does not admit a K\"ahler-Einstein metric then there is some $\beta\leq 1$ and a test configuration ${\cal X}$ with $\Fut_{\beta}({\cal X})=0$.

The problem of deforming the cone angle is well-understood \cite{kn:sd2} (there is a technical point concerning the holomorphic automorphism group which is clarified by Song and Wang in \cite{kn:SW}). So the problem becomes to show that if we have an increasing sequence $\beta_{i}\rightarrow \beta_{\infty}$, with $\beta_{\infty}\leq 1$ such that there are K\"ahler-Einstein metrics $\omega_{i}$ with these angles then {\it either} there is a K\"ahler-Einstein metric with cone angle $\beta_{\infty}$ {\it or} there is a test configuration ${\cal X}$ with $\Fut_{\beta_{\infty}}({\cal X})=0$ (with the same remark as above about the small complication when $X$ has automorphisms).

This strategy can be regarded as a variant of the standard \lq\lq continuity method'', in which one perturbs the K\"ahler-Einstein equation using a positive (1,1) form. In the cone-singularity case the current defined by the divisor takes the place of this (1,1) form. The advantages of the cone-singularity method are, first, that
  the solutions, for all $\beta$, have am intrinsic differential-geometric meaning.
  This allows us to extend many of the results proved for smooth K\"ahler-Einstein metrics to the singular case, where the corresponding extension to the continuity method is not known.  The extensions we outline below focus on holomorphic sections and the H\"ormander technique. Another important fact is that the $L^{2}$ norm of the curvature is a topological invariant \cite{kn:SW}. This was an important motivation for us in developing the approach (along the lines of \cite{kn:CD}) although in the end we do not use these ideas in our proof of Theorem 1.  A second advantage of the cone-singularity method is that there is a straightforward connection with algebraic geometry, as in the definition of $\Fut_{\beta}({\cal X})$ above.

\subsection{Sketch of more technical results}

Since there is a complete existence theory for the case of non-positive Ricci curvature, we can suppose that the Ricci curvature of $\omega_{i}$ is positive, say ${\rm Ric}(\omega_{i})= c_{i}\omega_{i}$ where $c_{i}\geq c>0$. Our first result is that we can approximate these metrics arbitrarily closely,  in the Gromov-Hausdorff sense, by smooth K\"ahler metrics $\omega'_{i}$ with ${\rm Ric}(\omega'_{i})\geq c'_{i}$ where $c'_{i}$ is as close as we like to $c_{i}$.
This means that we can immediately transfer the general Cheeger-Colding convergence theory for metrics of positive Ricci curvature  to the singular metrics. Thus, perhaps taking a subsequence, we can suppose that $(X_{i},\omega_{i})$ have a Gromov-Hausdorff limit $Z$(which is, initially,  just a length space).

The next step is to adapt the results and techniques of \cite{kn:DS} to show that $Z$ carries a natural algebraic structure.(As mentioned in \cite{kn:DS}, the application to metrics with cone singularities was in fact initiated before the more straightforward case treated there.)  The discussion now divides into two cases, when $\beta_{\infty}<1$ or when $\beta_{\infty}=1$. There are different difficulties in the two cases. An important feature in both cases is that the \lq\lq H\"ormander technique'' is used both \lq\lq globally'', to construct the projective embedding, and \lq\lq locally'', to study the local structure of the limit.  

In the case when $\beta_{\infty}<1$ it is clear that the \lq\lq regular set''
in $Z$ is open, the limiting metric is smooth there and the convergence is in $C^{\infty}$ there. The same applies to tangent cones to $Z$. We say such a tangent cone $C(Y)$ is \lq\lq good'' if the following holds. For any $\eta>0$
there is a function on $Y$ supported in the $\eta$-neighbourhood of the singular set $S(Y)$, equal to $1$ on some neighbourhood of $S(Y)$ and with the $L^{2}$ norm of its derivative less than $\eta$. One main technical result is that in fact all tangent cones are good. Given this, the arguments of \cite{kn:DS} extend almost word-for-word. The proof of this \lq\lq good tangent cones'' property involves a somewhat complicated argument, revolving around the  local density of the singular set $D$ in $(X,\omega_{i})$ (i.e. the volumes of the intersection of the divisor with small metric balls).

In the case when $\beta_{\infty}=1$ the main problem is to show that the regular set is open and the limiting metric is smooth there. We have two arguments for this: one involving again a study of the local densities of the singular set and the other using Ricci flow. 

The upshot in either case is that we are able to show that $Z$ is homeomorphic to a normal projective variety $W$ with log terminal singularities, just as in \cite{kn:DS}. Moreover the Gromov-Hausdorff convergence can be mirrored by algebro-geometric convergence in the standard sense. That is, we can find some fixed large $r$ such that if we write $N= {\rm dim} H^{0}(X; K_{X}^{-r}) -1$ then there is a sequence of embeddings $\iota_{i}:X\rightarrow \bC\bP^{N}$ each realised by the complete linear system $\vert - r K_{X}\vert$ such that
$\iota_{i} X$ converges to $W$ as projective varieties. Moreover, we get a limiting divisor $\Delta$ in $W$. Now there are three possibilities
\begin{enumerate}
\item The pair $(W,\Delta)$ is isomorphic  to $(X,D)$.
\item The pair $(W,\Delta)$ is not isomorphic to $(X,D)$ but $W$ is isomorphic to $X$.
\item  $W$ is not isomorphic to $X$.
\end{enumerate}
The second case further divides into (2a): when $\beta_{\infty}=1$ and (2b): when $\beta_{\infty}<1$.

To complete the proof we need to show that in cases (1) and (2a) there is a K\"ahler-Einstein metric on $(X,D)$ with cone angle $\beta_{\infty}$ and in cases (2b) and (3) that  the pair $(W,\Delta)$ is the central fibre $(X_{0}, D_{0})$ of a test configuration ${\cal X}$ with $\Fut_{\beta_{\infty}}({\cal X})=0$. Then in case (3) ${\cal X}$ is non-trivial and $\Fut({\cal X})\leq 0$, contradicting stability. In case (2b) ${\cal X}$ is trivial but, since $\beta_{\infty}<1$, we get $\Fut({\cal X})<0$, again contradicting stability.  

Case (1) can be handled in two ways. One is to show that the K\"ahler potentials, relative to a fixed reference metric, of $\omega_{i}$  are bounded, then to use the $C^{2}$ estimate of \cite{kn:JMR}. The other operates with the convergence theory to show that when the  complex structure of the limiting pair $(W,\Delta)$ is smooth the limiting metric has standard cone singularities. Case (2a) is similar (since $\beta_{\infty}=1$ there is no singularity in the metric).

 Case (3), which is the crucial issue,  involves two main difficulties. (Case (2b) is covered by the same discussion.) On the algebro-geometric side, we can think of $\iota_{i}(X,D)$ as a sequence in a fixed $PGL(N+1)$-orbit  in the appropriate Chow variety and $(W,\Delta)$ as a point in the closure of that orbit. Saying that $(W,\Delta)$ is the central fibre of a test configuration is the same as saying that $(W,\Delta)$ lies in the closure of a $\bC^{*}$-orbit for a suitable $1$-parameter subgroup in $PGL(N+1)$. There is no general reason why this should be true---see the discussion in \cite{kn:sd1} for example.  However it is true, by the Luna Slice theorem and the Hilbert-Mumford theorem applied to a slice, if we know that the automorphism group of $(W,\Delta)$ is {\it reductive}. Thus we need to prove
\begin{itemize}
\item $Aut(W,\Delta)$ is reductive.
\item The Futaki invariant $\Fut_{\beta_{\infty}}$ vanishes.
\end{itemize}
The proofs of these facts use  in an essential way recent developments in the theory of K\"ahler-Einstein metrics on singular spaces
 Although the detailed structure of the limiting metric $\omega_{\infty}$ on $Z$ seems quite difficult to study, we are able to show that it fits into the class of metrics considered in \cite{kn:BBEGZ}, \cite{kn:Bern}, \cite{kn:Ber2} and is a weak solution of the K\"ahler-Einstein equation (perturbed, of course, by a contribution from $(\Delta, \beta_{\infty})$). It can be characterised as a critial point of the perturbed Ding functional.  Then the  uniqueness result of Berndtsson \cite{kn:Bern}, as extended in \cite{kn:BBEGZ},   can be used to show that the automorphism group is reductive. This is a variant of the standard Matsushima Theorem, which asserts that the automorphism group of a manifold with a smooth K\"ahler-Einstein metric is reductive; the new feature being that the proof operates with the Lie groups rather than their Lie algebras. In a similar vein, the vanishing of the Futaki invariant follows from the  recent results of Berman \cite{kn:Ber2}, related to earlier results of Ding and Tian \cite{kn:DT}.

\end{document}